%
%
%

\input amstex
\documentstyle{amsppt}

\hsize=4.75in
\vsize=8in

\font\tencirc=lcircle10
\newbox\frogup
\newbox\frogdown
\long\def\hookdownarrow{%
\setbox\frogup=\hbox{\tencirc \char3 \char0}%
\setbox\frogdown=\hbox to-0.25\wd\frogup{\hss \hbox to 0pt{\hss 
$\displaystyle \downarrow $\hss }\kern1.9pt}%
\vrule height0pt width0pt depth1.5\ht\frogdown%
\raise\dp\frogup\box\frogup\lower\ht\frogdown\box\frogdown}
\font\ninerm=cmr9
\define\boldone{\hbox{\rm{\ninerm1\kern-3.3pt}1}}
\define\bZ{{\Bbb Z}}
\define\bz{{\Bbb Z}}

\define\bT{{\Bbb T}}
\define\bR{{\Bbb R}}

\define\bC{{\Bbb C}}

\define\cV{{\Cal V}}
\define\cT{{\Cal T}}
\define\cO{{\Cal O}}
\define\cH{{\Cal H}}
\define\cK{{\Cal K}}
\define\cF{{\Cal F}}
\define\cB{{\Cal B}}
\define\cD{{\Cal D}}
\define\relaxref{}
\define\relaxeqb{}
\define\relaxeqe{}
\define\relaxeqbnono{}
\define\relaxeqenono{}
\define\em{\it}
\define\mbox{\text}
\define\spaceddot{{\,\cdot\,}}
\newdimen\displayboxwidth
\displayboxwidth\hsize
\advance\displayboxwidth by-1.5in

\rightheadtext {Multiresolution Wavelet Theory and the Cuntz Algebra}

\topmatter
\title
A Connection between Multiresolution Wavelet Theory of
Scale $N$ and Representations of the Cuntz Algebra $\cO_N$
\endtitle
\author
Ola Bratteli and Palle E.T. Jorgensen
\endauthor
\affil
Mathematics Institute, 
University of Oslo,\\
PB 1053 --- Blindern, 
N-0316 Oslo, Norway\\
{\tt bratteli\@math.uio.no}
\\ \\
Department of Mathematics, 
University of Iowa,\\
14 MacLean Hall, 
Iowa City, IA 52242-1419, USA\\
{\tt jorgen\@math.uiowa.edu}
\endaffil
\thanks
Supported by the US NSF and the Norway Research Council.
\endthanks
\endtopmatter

\document

\head
1. Introduction
\endhead

In this paper we will give a short survey of a connection between the
theory of wavelets in $L^2(\bR)$ and certain representations of the
Cuntz algebra on $L^2(\bT)$.  This connection was first pointed out in
\cite{12} and has been developed further in \cite{5} and \cite{6}, and
these references contain complete proofs.  Basic reference for wavelet
theory is \cite{9} and for the Cuntz algebra
\cite{7}.  Let us emphasize at the outset that this is a field with more
questions than answers, and even quite fundamental questions are wide
open.  For example, it is hard to pinpoint abstractly which
representations of $\cO_N$ are obtained (although they can be
``written down'' as we will see), and it is unclear how various
equivalence relations between wavelets that one may envisage (same
father function up to scaling and translation, etc.) are reflected in
equivalence relations between representations (unitary equivalence,
quasi-equivalence, etc.).  The decomposition theory of the
representations has not been obtained in general, although it has been
worked out in great detail for related representations in \cite{2},
\cite{4}, \cite{5}, \cite{6}, and \cite{8}.

\head
2. From wavelets to representations
\endhead

Since wavelet theory of scale $N$ seems non-standard in the literature 
(but see \cite{11}), we will give
it a rundown here (see \cite{6, Section 10} for proofs).  Define
scaling by $N$ on
$L^2(\bR)$ by the unitary operator
\relaxeqb$$
(U\xi)(x)=N^{-\frac 12}\xi (N^{-1}x)
\text{\quad for } \xi\in L^2(\bR), x \in\bR.\tag{2.1} 
\relaxeqe$$
Let the {\em father function} be a unit vector $\varphi\in L^2(\bR)$,
and let $\cV_0$
be the closed linear span of the translates $T^k\varphi$, $k\in\bZ$, where
\relaxeqb$$
(T^k\xi)(x)=\xi(x-k)\tag{2.2}
\relaxeqe$$
is translation by $k$.  One assumes that $\varphi$ has the properties
\relaxeqb$$
\align
&\hbox to\displayboxwidth{$\displaystyle \{T^k\varphi\}_{k\in\bZ}$
is an orthonormal set in $L^2(\bR)$,\hss}\tag{2.3}
\\
&U\varphi \in\cV_0,\tag{2.4}
\\
&\bigwedge_{n\in\bZ} U^n \cV_0=\{0\},\tag{2.5a}
\\
&\bigvee_{n\in\bZ} U^n\cV_0 = L^2(\bR).\tag{2.5b}
\endalign
\relaxeqe$$
One example is the Haar father function:  $\varphi(x)=\chi_{[0,1]}(x)$.  By
(\relaxref{2.3}) we may define an isometry $\cF_\varphi :\cV_0\rightarrow
L^2(\bT):
\xi\rightarrow m$ as follows:  if
\relaxeqb$$
\xi(\spaceddot)=\sum_n b_n \varphi (\spaceddot -n)\tag{2.6}
\relaxeqe$$
then
\relaxeqb$$
m(t)=m(e^{-it})=\sum_n b_n e^{-int}\tag{2.7}
\relaxeqe$$
and we have the connection
\relaxeqb$$
\hat{\xi}(t)=m(t)\hat{\varphi}(t),\tag{2.8}
\relaxeqe$$
where $\xi \rightarrow \hat{\xi}$ is the Fourier transform, normalized so that
$\left\| \xi \right\|_{2}=\left\| \smash{\hat{\xi}}\vphantom{\xi}\right\|_{2}$.
In particular, if $\xi\in\cV_{-1}=U^{-1}\cV_0$, then $U\xi\in\cV_0$, and then
\relaxeqb$$
m_\xi =\cF_\varphi (U\xi)\in L^2(\bT)\ ,\tag{2.9}
\relaxeqe$$
and
\relaxeqb$$
\sqrt{N}\hat{\xi} (Nt)=m_\xi (t)\hat{\varphi}(t)\ .\tag{2.10}
\relaxeqe$$
In particular, using (2.4), we define
\relaxeqb$$
m_0(t)=m_\varphi (t)\ .\tag{2.11}
\relaxeqe$$
Now the condition (\relaxref{2.3}) is equivalent to
\relaxeqb$$
\text{\rm PER}
(|\hat{\varphi}|^2)(t)=\sum_k \left|\hat{\varphi}(t+2\pi k)\right|^2 =
(2\pi)^{-1},\tag{2.12}
\relaxeqe$$
and this implies
\relaxeqb$$
\sum^{N-1}_{k=0} \left|m_0 (t+2\pi {k/N})\right|^2=N\ . \tag{2.13}
\relaxeqe$$
If $\xi$, $\eta\in U^{-1}\cV_0$, then $\xi$ and $T^k\eta$ are orthogonal
for all
$k\in\bZ$ if and only if
\relaxeqb$$
\sum^{N-1}_{k=0} \bar{m}_\xi (t+2\pi{k/N}) m_\eta (t+2\pi {k/N})=0\tag{2.14}
\relaxeqe$$
for almost all $t\in\bR$, and $\{\xi (\spaceddot -k)\}_{k\in\bZ}$
is an orthonormal set
if and only if
\relaxeqb$$
\sum^{N-1}_{k=0}\left| m_\xi (t+2\pi {k/N})\right|^2=N. \tag{2.15}
\relaxeqe$$
With $m_0$ already given, we now choose $m_1,\dots,m_{N-1}$ in
$L^2(\bT)$ such that
\relaxeqb$$
\sum^{N-1}_{k=0} \bar{m}_i(t+2\pi{k/N}) m_j(t+2\pi {k/N}) =
\delta_{ij} N\tag{2.16}
\relaxeqe$$
for all $t\in\bR$, $i,j=0,\dots,N-1$.  If we define
$\psi_1,\dots,\psi_{N-1}$ by
\relaxeqb$$
\sqrt{N}\hat{\psi}_i (Nt) = m_i(t)\hat{\varphi} (t)\tag{2.17}
\relaxeqe$$
for $t\in\bR$, $i=1,\dots,N-1$, it follows that $\{ T^k\psi_i\}_{k\in\bZ,
i\in\{1,\dots,N-1\}}$ form an orthonormal basis for $\cV_{-1}\cap 
\cV^\perp_0$, and
hence $\{U^n T^k\psi_i\}_{k\in\bZ,i\in\{1,\dots,N-1\}}$ form an
orthonormal basis
for $L^2(\bR)$.  The functions $\psi_1,\dots,\psi_{N-1}$ are called {\em mother
functions}.  They are not unique, but depend on the choice of the functions
$m_1,\dots,m_{N-1}$ satisfying (\relaxref{2.16}).

If $\rho=\rho_N=e^{\frac {2\pi i}N}$, the condition (\relaxref{2.16})
translates into the
requirement that the $N\times N$ matrix
\relaxeqb$$
{\frac 1{\sqrt{N}}}
{\pmatrix
           m_0(z) & m_0(\rho z) & \dots & m_0(\rho^{N-1} z) \\
           m_1(z) & m_1(\rho z) & \dots & m_1(\rho^{N-1} z) \\
           \vdots & \vdots & \ddots & \vdots  \\
           m_{N-1}(z) & m_{N-1}(\rho z) & \dots & m_{N-1}(\rho^{N-1} z)
       \endpmatrix}\tag{2.18}
\relaxeqe$$
is unitary for almost all $z\in\bT$.

Now, this is again equivalent to saying that the operators $S_i$ defined on
$L^2(\bT)$ by
\relaxeqb$$
(S_i\xi)(z) = m_i(z)\xi (z^N)\tag{2.19}
\relaxeqe$$
for $\xi\in L^2(\bT)$, $z\in\bT$, $i=0,1,\dots,N-1$ satisfy the relations
\relaxeqb$$
\align
S^*_j S^{}_i&=\delta^{}_{ij}\boldone\tag{2.20} \\
\sum^{N-1}_{i=0} S^{}_i S^*_i &=\boldone,\tag{2.21}
\endalign
\relaxeqe$$
which are exactly the Cuntz relations.  This defines the map from the $N$-scale
multiresolution wavelet $\{\varphi, \psi_1,\dots,\psi_{N-1}\}$ into
representations
of $\cO_N$.

\head
3. From representations to wavelets
\endhead

When, conversely, does a representation of the Cuntz algebra $\cO_N$
give rise to a multiresolution wavelet
$\{\varphi,\psi_1,\dots,\psi_{N-1}\}$ such that one can recover the
representation again from the wavelet by the construction in Section
2?  A minimal requirement is that the representation acts on
$L^2(\bT)$ by formula (\relaxref{2.19}), and then unitarity of
(\relaxref{2.18}) is assured from the Cuntz relations.  For any
representation of $\cO_N$ on a Hilbert space $\cH$ we may define an
associated endomorphism $\sigma$ of $\cB(\cH)$ by $\sigma(\spaceddot)
= \sum^{N- 1}_{i=0} S^{}_i\spaceddot S^*_i$ (see, e.g., \cite{2}).
When $\cH=L^2(\bT)$, and $S_i$ is given by (\relaxref{2.19}), a simple
computation, using unitarity of (\relaxref{2.18}), shows that
\relaxeqb$$
\sigma(M_f)=M_{\bar{\sigma}(f)}\tag{3.1}
\relaxeqe$$
for all $f\in L^\infty(\bT)$, where $M_f$ is the operator of multiplication
by $f$ on
$L^2(\bT)$, and $\bar{\sigma}(f)(z)=f(z^N)$.  Conversely,
\proclaim{Proposition 3.1}
\rom{(\cite{6, Proposition 1.1})} If $S_0,\dots,S_{N-1}$ is a representation
of $\cO_N$
on $L^2(\bT)$ and
\relaxeqb$$
\sum^{N-1}_{i=0} S^{}_i M^{}_f S^*_i = M^{}_{\bar{\sigma}(f)}\tag{3.2}
\relaxeqe$$
for all $f\in L^\infty (\bT)$, then $S_i$ has the form \rom{(\relaxref{2.19})}
with
\relaxeqb$$
m_i=S_i \boldone.\tag{3.3}
\relaxeqe$$
\endproclaim
\demo{Proof}
If $f\in L^\infty(\bT)\subset L^2(\bT)$ then
\relaxeqbnono$$
\align
 M^{}_{\bar{\sigma}(f)} S^{}_j & = \sum_i S^{}_i M^{}_f S^*_i 
S^{}_j\tag{3.4} \\
& = S_j M_f,
\endalign
\relaxeqenono$$
and applying this to $\boldone$ and using (\relaxref{3.3}) we have
\relaxeqb$$
f(z^N)m_j(z) = (S_jf)(z).\tag{3.5}
\relaxeqe$$
As $L^\infty(\bT)$ is dense in $L^2(\bT)$, (\relaxref{2.19}) follows.
\enddemo

In order that the representation shall give rise to wavelets, it is
not sufficient that it have the form (\relaxref{2.19}), however, as we
will discuss further in the next section.  Let us for the moment
assume that the representation comes from a wavelet satisfying the
slight regularity condition that $\hat{\varphi}(t)$ is continuous near
$t=0$.  Then condition (\relaxref{2.5a}) implies $\hat{\varphi}(0)\ne
0$ (see
\cite{9, Remark 3 after Proposition 5.3.2}).  It follows from
(\relaxref{2.11}) and
(\relaxref{2.10}) that

\relaxeqb$$
\sqrt{N} \hat{\varphi} (Nt) = m_0 (t)\hat{\varphi} (t),\tag{3.6}
\relaxeqe$$
and hence $m_0(t)$ is continuous near $t=0$ and $m_0(0)=\sqrt{N}$.  
Thus it follows
from (\relaxref{2.13}) that
\relaxeqb$$
m_0(2\pi {k/N}) = 0\tag{3.7}
\relaxeqe$$
for $k=1,\dots,N-1$, and combining this with (\relaxref{3.6})
and using a recursive
argument we deduce that
\relaxeqb$$
\hat{\varphi}(2\pi k) = 0\tag{3.8}
\relaxeqe$$
for all $k\in\bZ\setminus \{0\}$.  It now follows from (\relaxref{2.12}) that
\relaxeqb$$
|\hat{\varphi} (0)| = (2\pi)^{-{\frac 12}},\tag{3.9}
\relaxeqe$$
and by changing $\hat{\varphi}$ by an irrelevant phase factor we may assume
$\hat{\varphi}(0)=(2\pi)^{-{\frac 12}}$.  But an iteration of
(\relaxref{3.6}) gives
\relaxeqb$$
\hat{\varphi}(t) = \prod^n_{k=1} (N^{-{\frac 12}} m_0 (tN^{-k}))
\hat{\varphi} (tN^{-
n}),\tag{3.10}
\relaxeqe$$
and as $\lim_{n\rightarrow\infty} \hat{\varphi} (tN^{-n}) = \hat{\varphi} (0)
= (2\pi)^{-{\frac 12}}$ we deduce
\relaxeqb$$
\hat{\varphi}(t) = (2\pi)^{-{\frac 12}} \prod^{\infty}_{k=1}
(N^{-{\frac 12}} m_0(tN^{-
k})).\tag{3.11}
\relaxeqe$$
Under even stronger regularity properties on $\varphi$, for example
that $m_0$ is Lipschitz continuous near $0$, the expansion
(\relaxref{3.11}) converges absolutely and uniformly on compact sets.

If we view functions $\xi$ in $L^2(\bT)$ as $2\pi$-periodic functions on
$\bR$, it
follows from (\relaxref{2.19}) that
\relaxeqb$$
(S^n_0 \xi)(t) = \prod^{n-1}_{k=0} m_0 (N^k t)\xi (N^n t),\tag{3.12}
\relaxeqe$$
and if $E:L^2(\bT)\rightarrow L^2 (\bR)$ is the embedding determined by
\relaxeqb$$
 (E\xi)(t) =\left\{ \aligned
\xi(t) & \quad\mbox{if $-\pi\leq t\leq\pi$} \\
     0 & \quad\mbox{otherwise}
\endaligned\right. \tag{3.13}
\relaxeqe$$
then
\relaxeqb$$
(U^n E S^n_0\xi)(t) = \chi_{[-\pi,\pi]}(tN^{-n}) \prod^n_{k=1} (N^{-{\frac 12}}
m_0(N^{-k}t))\xi (t).\tag{3.14}
\relaxeqe$$
Thus it follows from (\relaxref{3.11}) that
\relaxeqb$$
\lim_{n\rightarrow\infty} U^n_{} E S^n_0 \xi = (2\pi)^{\frac 12} \hat{\varphi}
\xi,\tag{3.15}
\relaxeqe$$
where the convergence is uniform on compact subsets of $\bR$ if
$\xi \in L^\infty
(\bT)\subset L^2 (\bT)$.  In a similar way, using (\relaxref{2.17})
and iteration, one
deduces
\relaxeqb$$
\lim_{n\rightarrow\infty} U_{}^n E S^{n-1}_0 S^{}_i \xi = (2\pi)^{\frac 12}
\hat{\psi}^{}_i
\xi.\tag{3.16}
\relaxeqe$$
Thus, the formulae (\relaxref{3.15}) and (\relaxref{3.16}) allow us to recover the wavelet
system $\{\varphi, \psi_1,\dots,\psi_{N-1}\}$ from the representation.

\head
4. Which representations may occur?
\endhead
Are there other criteria than those in Proposition 3.1 ensuring that a
representation of $\cO_N$ has the form (\relaxref{2.19})?  
A necessary condition can be
formulated in terms of the Wold decomposition of the isometries $S_i$.  
In general,
if $S$ is an isometry, define a decreasing sequence of projections by
\relaxeqb$$
E_k = S^k S^{*k}\tag{4.1}
\relaxeqe$$
and let $P_U$ be the limit projection
\relaxeqb$$
P_U = \mathop{\text{\rm s-lim}}_{k\rightarrow\infty} E_k.\tag{4.2}
\relaxeqe$$
Then $SP_U=P_US$, $P_US$ is a unitary operator on $P_U\cH$, and
$(1-P_U)S$ is a shift
on $(1-P_U)\cH$, i.e.,
\relaxeqb$$
\bigcap_n S^n (1-P_U)\cH = \{0\}.\tag{4.3}
\relaxeqe$$
The Wold decomposition is
\relaxeqb$$
S=SP_U\oplus S(1-P_U).\tag{4.4}
\relaxeqe$$

\proclaim{Theorem 4.1}
Let $S$ be an operator on $L^2(\bT)$ of the form
\relaxeqb$$
(S\xi)(z)=m(z)\xi (z^N)\tag{4.5}
\relaxeqe$$
and assume that $S$ is an isometry, i.e.,
\relaxeqb$$
\sum^{N-1}_{k=0} |m(\rho^k z)|^2=N\tag{4.6}
\relaxeqe$$
for almost all $z\in\bT$, where $\rho=e^{\frac {2\pi i}N}$.  It
follows that the projection $P_U$ corresponding to the unitary part of
the Wold decomposition of $S$ is one- or zero-dimensional.
Furthermore, it is one-dimensional if and only if both conditions
\rom{(\relaxref{4.7})} and \rom{(\relaxref{4.8})} are satisfied.
\relaxeqb$$
\align
&|m(z)|=1\quad\text{for almost all }z\in\bT.\tag{4.7} \\
&\text{There exists a measurable function }\xi : \bT\rightarrow\bT\tag{4.8} \\
&\text{and a }\lambda\in\bT\text{ such that} \\
&\hbox to\displayboxwidth{\hss$\displaystyle m(z)\xi (z^N)=\lambda\xi(z)$\hss}
\\
&\text{for almost all }z\in\bT.
\endalign
\relaxeqe$$
In this case the range of the projection $P_U$ is
$\bC\xi$.
\endproclaim

\demo{Proof}
See \cite{6, Theorem 3.1}.  This paper also contains more general versions of
Theorem 4.1.
\enddemo
Now, combining (\relaxref{4.7}) with (\relaxref{2.17}) and using the
ergodicity of $z\mapsto z^N$ one deduces:

\proclaim{Corollary 4.2}
The operators $S_0,\dots,S_{N-1}$ in the representation of $\cO_N$
defined by a wavelet $\{\varphi,\psi_1,\dots,\psi_{N-1}\}$ have zero
unitary part in the Wold decomposition, i.e., they are all shifts.
\endproclaim

\demo{Proof}
Lee Lemma 9.3 in \cite{6}.
\enddemo
Corollary 4.2 gives a rather severe restriction on the representations
that can be defined by wavelets.  In the same way as a single shift is
always isomorphic to a multiple of the shift given by multiplication
by $z$ on the Hardy space $\cH_+(L^2(\bT))$, one may use the shift
property of $S_i$ (actually it suffices that $S_0$ is a shift for the
following construction) to realize the representation $\{S_0,
\dots,S_{N-1}\}$ of $\cO_N$ on $\cK=L^2(\bT)$ on the Hilbert space.
\relaxeqb$$
\cH_+ \left(\bigoplus^{N-1}_{j=1}\cK\right) = \bigoplus^\infty_{n=1}
\left(\bigoplus^{N-
1}_{j=1}\cK\right),\tag{4.9}
\relaxeqe$$
where we view the elements $(\xi_n)\in \cH_+(\cK)$ as the 
$\bC^{N-1}\otimes \cK$-valued functions
\relaxeqb$$
\xi(z)=\sum^\infty_{n=1}\xi_n z^n\tag{4.10}
\relaxeqe$$
on $\bT$, such that $S_0$ is represented by the operator $M_z=$
multiplication by $z$.  To this end we define a unitary operator $V :
\cH_+ \left(\bigoplus^{N-1}_{j=1}
\cK\right)\rightarrow\cK$ by
\relaxeqb$$
V\left(\sum^\infty_{k=1} \left(\bigoplus^{N-1}_{j=1}
\psi^{(j)}_k\right)z^k_{}\right) = \sum^\infty_{k=1}
\sum^{N-1}_{j=1} S_0^{k-1} S^{}_j \psi^{(j)}_k\tag{4.11}
\relaxeqe$$
The Cuntz relation together with $\mathop{\text{\rm
s-lim}}_{k\rightarrow\infty} S^k_0 S^{*k}_0=0$ ensures that $V$ is
unitary (see Lemma 6.1 in \cite{6}), and if $S^+_i=V^* _{}S^{}_iV$ for
$i=0,\dots,N-1$, one verifies that $$S^+_0 = M_z =\text{
multiplication by }z$$ and
\relaxeqbnono$$
S^+_i\psi = z\left(\bigoplus^{i-1}_{j=1} 0\oplus V \psi \oplus 
\left(\bigoplus^{N-1}_{j=i+1} 0\right)\right)\text{\quad for }i=1,\dots,N-1.
\relaxeqenono$$
Now, the space $\cH_+
\left(\bigoplus^{N-1}_{j=1}\cK\right)=\bC^{N-1}_{}\otimes \cK\otimes
H^2_+ (\bT)$, where $H^2_+(\bT)$ consists of the functions $\xi\in
L^2(\bT)$ with a Fourier expansion of the form $\sum^\infty_{n=1}
a_nz^n$, has an obvious embedding in $\bC^{N-1}\otimes\cK\otimes
L^2(\bT)$.  If the representation comes from a wavelet
$\{\varphi,\psi_1,\dots,\psi_N\}$, so that $\cK=L^2(\bT)$, one may
define a unitary map
\relaxeqbnono$$
J:L^2(\hat{\bR})\rightarrow\bC^{N-1}\otimes\cK\otimes L^2(\bT)
\relaxeqenono$$
by the requirement
\relaxeqbnono$$
J(U^n T^k\psi_m)(e^{-it},z)=e_m\otimes e^{-ikt}\otimes z^n
\relaxeqenono$$
for $n,k\in\bZ$, $m=1,\dots,N-1$, where $\{e_m\}^{N-1}_{m=1}$ is the
standard basis in $\bC^{N-1}$.  One can then establish that the
diagram

\relaxeqbnono$$
\matrix
\cV_{0} & 
\mkern18mu \botsmash{\raise0.5ex\hbox to0pt{\hss $\displaystyle 
\overset{\cF_\varphi }\to{\longrightarrow }$\hss }\lower0.5ex\hbox 
to0pt{\hss $\displaystyle \underset{\cF_{\varphi }^{-1}}\to{\longleftarrow 
}$\hss }}\mkern18mu %
& \cK=L^{2}\left( \bT\right) & 
\mkern18mu \botsmash{\raise0.5ex\hbox to0pt{\hss $\displaystyle 
\overset{V^*}\to{\longrightarrow }$\hss }\lower0.5ex\hbox to0pt{\hss 
$\displaystyle \underset{\ }\to{\longleftarrow }$\hss }}\mkern18mu %
& \bC^{N-1}\otimes\cK\otimes
H_{+}^{2}\left( \bT\right) \\ 
\hookdownarrow %
&  & 
\hbox to0pt{\hss \vrule height6pt width0.26pt depth0pt \hss 
}\setbox\frogdown=\hbox to 0pt{\hss $\displaystyle \downarrow $\hss 
}\lower\ht\frogdown\box\frogdown\lower0.5ex\hbox to0pt{$\scriptstyle 
\mkern6mu M_{\hat{\varphi }}$\hss }%
\vrule height16pt width0pt depth20pt &  & 
\hookdownarrow %
\\ 
L^{2}\left( \bR\right) & 
\mkern18mu \topsmash{\raise0.5ex\hbox to0pt{\hss $\displaystyle 
\overset{\cF}\to{\longrightarrow }$\hss }\lower0.5ex\hbox to0pt{\hss 
$\displaystyle \underset{\cF^{-1}}\to{\longleftarrow }$\hss }}\mkern18mu %
& L^{2}\left( 
\smash{\hat{\bR}}\vphantom{\bR}%
\right) & 
\mkern18mu \topsmash{\raise0.5ex\hbox to0pt{\hss $\displaystyle 
\overset{J}\to{\longrightarrow }$\hss }\lower0.5ex\hbox to0pt{\hss 
$\displaystyle \underset{J^{*}}\to{\longleftarrow }$\hss }}\mkern18mu %
& \bC^{N-1}\otimes\cK\otimes L^{2}\left( \bT\right)
\endmatrix
\relaxeqenono$$
is commutative, where $\cF_\varphi$ is defined prior to
(\relaxref{2.6}) and $\cF$ is Fourier transform, so that the left
rectangle in the diagram is commutative by (\relaxref{2.8}).  Note
that the scaling operator $U$ on $L^2(\bR)$ is transformed into the
operator of multiplication by $z$ on the space
$\bC^{N-1}\otimes\cK\otimes L^2(\bT)$ by conjugation by the unitary
$J\cF$.

\head
5. Some other representations
\endhead
In \cite{2}, \cite{4}, \cite{3}, \cite{12}, \cite{5}, \cite{8}, and
\cite{6}, other representations of $\cO_N$, many of which are of the
form (\relaxref{2.19}) with unitarity of (\relaxref{2.18}), have been
considered.  For example, it is proved in \cite{5} that if the $m_i$
are monomials in $z$, that is,
\relaxeqb$$
m_i(z)=z^{d_i},\tag{5.1}
\relaxeqe$$
where $d_0,d_1,\dots,d_{N-1}$ are integers mutually incongruent modulo
$N$, then the representation has a discrete decomposition into
mutually disjoint irreducible subrepresentations, and the
corresponding irreducible subspaces are spanned by the monomials they
contain.  Thus the explicit form of the decomposition depends on
specific number-theoretic properties of the sequence
$d_0,\dots,d_{N-1}$, and has been worked out in full detail for $N=2$
in \cite{5}.  However, by Corollary 4.2, none of these representations
comes from a wavelet.  (If they are modified by the canonical action
of $U(N)$ on $\cO_N$, some of them do come from wavelets, for example
the Haar wavelet.)

Let $s_0,\dots,s_{N-1}$ denote the generators of $\cO_N$, satisfying
the relations (2.20) and (2.21).  Let $\text{\rm UHF}_N$ denote the
closure of the linear span of elements of the form $s^{}_{i_1}\cdots
s^{}_{i_n} s^*_{j_n}\cdots s^*_{j_1}$, and $\cD_N$ the closure of the
linear span of elements of the form $s^{}_{i_1}\cdots s^{}_{i_n}
s^*_{i_n}\cdots s^*_{i_1}$.

It is well known that $\text{\rm UHF}_N$ is a $\text{\rm UHF}$ algebra
of Glimm type $N^\infty$, and $\cD_N$ is a maximal abelian subalgebra
of $\text{\rm UHF}_N$, and of $\cO_N$.  If the monomial
representations mentioned in the previous paragraph are restricted to
$\text{\rm UHF}_N$, they have a similar decomposition theory to that
of the original representation, and there is a certain canonical
action of $\bz$ on the $\text{\rm UHF}_N$ irreducible components such
that the orbits of these components correspond to the $\cO_N$
irreducible components.

Let us consider another class of representations of the form (\relaxref{2.19}):

\proclaim{Proposition 5.1}
\rom{(This is Theorem {7.1} from \cite{6}.)}  Consider a representation
$\pi$ of $\cO_N$
on $L^2(\bT)$ of the form \rom{(\relaxref{2.19}):}
\relaxeqbnono$$
(S_i\xi)(z) = m_i(z)\xi(z^N),
\relaxeqenono$$
where the $m_i$ satisfy the unitarity condition \rom{(\relaxref{2.18}).} Let
$M_{L^\infty(\bT)}$ be the image of $L^\infty(\bT)$ acting as multiplication
operators on $L^2(\bT)$.  The following conditions are equivalent:
\relaxeqb$$
\align
&\hbox to\displayboxwidth{$\displaystyle \pi(\cD_N)''\subset 
M_{L^\infty(\bT)}$,\hss }\tag{5.2} \\
&\pi (\cD_N)'' = M_{L^\infty(\bT)},\tag{5.3} \\
&m_j(z)=\sqrt{N} \chi_{A_j}(z)u(z),\tag{5.4}
\endalign
\relaxeqe$$
where $u:\bT\rightarrow\bT$ is a measurable function, and
$A_0,A_1,\dots,A_{N-1}$ are $N$ measurable subsets of $\bT$ with the
property that if $\rho=\rho_N=e^{\frac {2\pi i}N}$, then for almost
all $z\in \bT$ the $N$ equidistant points $z,\rho z, \rho^2z,
\dots,\rho^{N-1} z$ lie with one in each of the sets
$A_0,\dots,A_{N-1}$ \rom{(}i.e., $A_0,\dots, A_{N-1}$ form a partition
of $\bT$ up to null-sets, and for each $k$ the sets $A_k, \rho
A_k,\dots,\rho^{N-1} A_k$ form a partition of $\bT$\rom{).}  Any $m_i$
of this form does indeed define a representation of $\cO_N$.
\endproclaim

To analyze these representations further, note that by a decoding on
$\bT$, i.e., a measure-preserving transformation of $\bT$, we may
assume that $A_k$ is the segment of $\bT$ between $\rho^k$ and
$\rho^{k+1}$, and then we put
\relaxeqb$$
\chi_k(z)=\chi_{A_k}(z).\tag{5.5}
\relaxeqe$$
Now let $S^{(j)}_k$, $j=1,2,\;k=0,\dots,N-1$ be two representations of
this kind, i.e., there exist measurable functions
$u^{(j)}:\bT\rightarrow \bT$ such that
\relaxeqb$$
(S^{(j)}_k\xi)(z)=\sqrt{N} \chi_k(z)u^{(j)}(z)\xi(z^N).\tag{5.6}
\relaxeqe$$
Then one can show that $T\in \cB(L^2(\bT))$ intertwines the two
representations, i.e.,
\relaxeqb$$
T\pi^{(1)}(x)=\pi^{(2)} (x)T\quad\text{for all }x\in\cO_N,\tag{5.7}
\relaxeqe$$
if and only if
\relaxeqb$$
\align
&T=M_f\text{ where }f\in L^\infty(\bT)\text{ is a function
satisfying}\tag{5.8} \\
&\hbox to\displayboxwidth{\hss$\displaystyle f(z)u^{(1)}(z) = 
u^{(2)}(z) f(z^N)$\hss} \\
&\text{for almost all }z\in\bT\text{ (i.e., the cocycles }u^{(1)}, 
u^{(2)}\text{ cobound} \\
&\text{with the coboundary }f\text{ for the action }z\mapsto z^N\text{).} 
\endalign
\relaxeqe$$
(See
\cite{6, Proposition 8.1}.)
Since the map $z\mapsto z^N$ is ergodic (w.r.t. Haar measure on
$\bT$), it follows that $f$ is unique up to a scalar multiple of a
function $\bT\rightarrow \bT$ if a nonzero $f$ exists at all.  In
particular, if $u^{(1)}=u^{(2)}$ then $f(z)=f(z^N)$ and $f$ is a
scalar multiple of $1$.  Thus (see \cite{6, Corollary 8.3}):

\proclaim{Corollary 5.2}
If the representations $\pi^{(j)}$, $j=1,2,$ are defined by 
\rom{(\relaxref{5.6}),} then
$\pi^{(j)}$ are irreducible, and the following conditions are equivalent.
\relaxeqb$$
\align
&\pi^{(1)}\text{ and }\pi^{(2)}\text{ are unitarily equivalent.}\tag{5.9} \\
&\text{The cocycles }u^{(1)}, 
u^{(2)}\text{ cobound, i.e., there exists a}\tag{5.10} \\
&\text{measurable function }\Delta :\bT\rightarrow\bT\text{ such that} \\
&\hbox to\displayboxwidth{\hss$\displaystyle \Delta (z) u^{(1)}(z) = 
u^{(2)} (z) \Delta(z^N)$\hss} \\
&\text{for almost all }z\in\bT.
\endalign
\relaxeqe$$
\endproclaim

Let us end this section by mentioning a completely different way of
describing states and representations of $\cO_N$ from \cite{4} and
\cite{3}.  If $\hat{\omega}$ is a state of $\cO_N$, $\pi$ the
associated representation on $\cH$ with cyclic vector $\Omega$ and
$S_i=\pi(s_i)$, let $\cK$ be the closed linear span of vectors of the
form $S^*_{i_1}\cdots S^*_{i_k}\Omega$, for
$k=0,1,\dots,\;i_j\in\{0,\dots,N-1\}$.  For example, if $\omega$ is a
Cuntz state (coherent state), $\cK$ is one-dimensional, and
conversely.  Let $P$ be the projection from $\cH$ onto $\cK$, and
define
\relaxeqbnono$$
V^*_i=PS^*_i P=S^*_iP\in \cB(\cH).
\relaxeqenono$$
The Cuntz relations (\relaxref{2.20})--(\relaxref{2.21}) imply that
\relaxeqbnono$$
\sum^{N-1}_{i=0} V^{}_i V^*_i = \boldone_\cK,
\relaxeqenono$$
and if
$\omega=\hat{\omega}|_{\cB(\cK)}$,
then
\relaxeqbnono$$
\hat{\omega}(s^{}_{i_1}\cdots s^{}_{i_n}s^*_{j_m}\cdots s^*_{j_1}) = 
\omega (V^{}_{i_1}\cdots
V^{}_{i_n} V^*_{j_m}\cdots V^*_{j_1})
\relaxeqenono$$
so $\hat{\omega}$ is completely determined by the pair $\omega,
\{V_i\}$.  Conversely:

\proclaim{Theorem 5.3}
\rom{(Popescu's reconstruction theorem \cite{13}, \cite{14}.)}  

If $V_0,\dots,V_{N-1}\in
\cB(\cK)$, where $\cK$ is a Hilbert space, $\Omega\in\cK$ is a unit
vector cyclic under
the polynomials in $V^*_0,\dots,V^*_{N-1}$, and $\sum^{N-1}_{i=0}
V^{}_iV^*_i=\boldone_\cK$, then there exists a state $\hat{\omega}$ on
$\cO_N$ such that
\relaxeqbnono$$
\hat{\omega} (s^{}_{i_1}\cdots s^{}_{i_n} s^*_{j_m}\cdots s^*_{j_1}) = \langle
V^*_{i_n}\cdots V^*_{i_1}\Omega|V^*_{j_m}\cdots V^*_{j_1}\Omega\rangle,
\relaxeqenono$$
that is, the $V_i$'s have a dilation to a representation of the $s_i$'s.
\endproclaim

The paper \cite{3} contains characterizations of pure states
$\hat{\omega}$, and states $\hat{\omega}$ with
$\pi_{\hat{\omega}}(\text{\rm UHF}_N)''= \cB(\cH_{\hat{\omega}})$, in
terms of ergodicity properties of the completely positive map
{$\underline{\sigma}$} $(\spaceddot)=\sum_k V^{}_k\spaceddot V^*_k$ of
$\cB(\cK)$.  For example, $\hat{\omega}$ is pure if and only if
$\underline{\sigma}$ is ergodic in the sense that
$\underline{\sigma}(X)=X$ implies $X\in \bC\boldone_{\cK}$, and this
is again true if and only if $\{V^{}_i,V^*_i\}$ acts irreducibly on
$\cK$ and the projection $P:\cH\rightarrow\cK$ is contained in
$\pi_{\hat{\omega}}(\cO_N)''$. Furthermore,
$\pi_{\hat{\omega}}(\text{\rm UHF}_N)''=\cB(\cH_{\hat{\omega}})$ if
and only if $\text{\rm Tail}(\underline{\sigma})=\bC\boldone_{\cK}$,
i.e., all w*-limit points of sequences of the form
$\underline{\sigma}^{n_k}(X_k)$, where $n_k\rightarrow\infty$ and
$X_k\in\cB(\cK)$ are uniformly bounded, are multiples of
$\boldone_{\cK}$.  Let us end by citing a proof of Reinhard Werner of
Theorem 5.3, which is substantially more direct than the original
proof in \cite{13}, \cite{14}: With the assumptions in Theorem 5.3, it
suffices by Stinespring's theorem to show that the map
$R:\cO_d\rightarrow\cB(\cK)$ defined by $R(\boldone)=\boldone$ and
\relaxeqbnono$$
R(s^{}_{i_1}\cdots s^{}_{i_n} s^*_{j_m}\cdots s_{j_1}) = 
V^{}_{i_1}\cdots V^{}_{i_n} V^*_{j_m}\cdots
V^*_{j_1}
\relaxeqenono$$
is completely positive.  Let $\cT_N$ be the Cuntz-Toeplitz algebra,
i.e., $\cT_N$ is the ${}^*$-algebra generated by $N$ isometries
$s_0,\dots,s_N$ with orthogonal ranges.  It is well known that $\cT_N$
is an extension of $\cO_N$ by the compact operators; see \cite{10}.
$\cT_N$ has a realization on the unrestricted Fock space $\hat{\cH} =
\bigoplus^\infty_{k=0} (\bC^N)^{\otimes k}$ by $\pi_0(s_i)\xi =
|e_i\rangle\otimes\xi$, where $\{ |e_i\rangle\}^{N-1}_{i=0}$ is the
standard basis of $\bC^N$.  Let $\lambda\in\bC$ with $|\lambda|<1$.
Define
\relaxeqbnono$$
W_\lambda : \cK \longrightarrow \hat{\cH}\otimes\cK
\relaxeqenono$$
by
\relaxeqbnono$$
W_\lambda\varphi = \sqrt{1-|\lambda|^2} \bigoplus^\infty_{k=0} 
\lambda^k \sum_{i_1\dots
i_k} |i_1\rangle\otimes\cdots\otimes |i_k\rangle \otimes 
V^*_{i_k}\cdots V^*_{i_1}
\varphi.
\relaxeqenono$$
One checks that $W_\lambda$ is an isometry, and
\relaxeqbnono$$
(\pi^{}_0(s^{}_i)^*\otimes\boldone^{}_\cK)W^{}_\lambda =\lambda 
W^{}_\lambda V^*_i.
\relaxeqenono$$
Define
\relaxeqbnono$$
\align
R_\lambda (s^{}_{i_1}\cdots s^{}_{i_n} s^*_{j_m}\cdots s^*_{j_1}) 
&:= W^*_\lambda (\pi^{}_0
(s^{}_{i_1}\cdots s^{}_{i_n} s^{*}_{j_m}\cdots s^{*}_{j_1})\otimes 
\boldone^{}_\cK)W^{}_\lambda  \\
&=
\bar{\lambda}^n_{} \lambda^m_{} V^{}_{i_1}\cdots V^{}_{i_n} 
V^*_{j_m}\cdots V^*_{j_1}.
\endalign
\relaxeqenono$$
It follows from this explicit Stinespring representation that $R_\lambda$ is
completely positive for each $\lambda$ with $|\lambda|<1$, and letting $\lambda
\rightarrow 1$ it follows that $R$ is completely positive as a map
from $\cT_N$ into $\cB(\cK)$.  It remains to show that $R$ defines a
map of the quotient $\cO_N$ of $\cT_N$, i.e., that $R$ annihilates the
ideal generated by the projection $\boldone - \sum^{N-1}_{i=0} s^{}_i
s^*_i$.  But this is easily checked from $\sum_i V^{}_i V^*_i =
\boldone^{}_\cK$.

\head
6. Further results and problems
\endhead
Consider $m_0,m_1\in L^\infty(\bT)$ given, and assume that the matrix
\relaxeqb$$
C(z)={\frac 1{\sqrt{2}}}\pmatrix
m_0(z) & m_1(z)\\
m_0(-z) & m_1(-z)\endpmatrix\tag{6.1}
\relaxeqe$$
is unitary for almost all $z\in\bT$.  Then consider the spectral
problem of finding $L^2(\bT)$-solutions $\varphi$ to
\relaxeqb$$
{\frac 1{\sqrt{2}}} (m_0(z)\varphi (z^2)+m_1(z)\varphi 
(-z^2))=\lambda\varphi (z)\
.\tag{6.2}
\relaxeqe$$
Recall that, with unitarity of (\relaxref{6.1}) assumed, the operators
\relaxeqb$$
(S_i\xi)(z)=m_i(z)\xi(z^2)\tag{6.3}
\relaxeqe$$
define a representation of $\cO_2$ acting on the Hilbert space
$L^2(\bT)$, and we are concerned, in \cite{5}, with the possible
decompositions of this class of representations.  For our analysis in
\cite{6}, we introduce a new {\em index},
\relaxeqbnono$$
\mbox{ind} (\pi):=\sum_\lambda \mbox{(the dimension of the space
of solutions to
(6.2))\ .}
\relaxeqenono$$
It is motivated in part by Arveson's index \cite{1}.  We denote by
$\pi$ the representation given by $\pi(s_i)=S_{i_1}$, where $S_i$ are
defined by (6.3).  As further motivation, note that, for any two
solutions $\varphi, \psi$ to (\relaxref{6.2}), the associated function
$\{\varphi,\psi\}$ defined by
\relaxeqbnono$$
z\mapsto \bar{\varphi}(z)\psi(z)+\bar{\varphi}(-z)\psi(-z)
\relaxeqenono$$
is necessarily constant on $\bT$ (a.e.).  We show that the index must
take on values as follows:
\relaxeqbnono$$
\mbox{ind} (\pi)=p\in \{0,1,2\}\ ,
\relaxeqenono$$
and then the representation $\pi$ contains $\rho_1\oplus\cdots\oplus
\rho_p$ where each representation $\rho_i$ is isomorphic to one given
by the Haar wavelet.  (It is understood that the sum is empty if
$p=0$.)

In a future paper, we plan to study and refine our new index, with a
view to picking up copies of isomorphism classes of wavelets other
than the Haar one.  Certainly the various families of wavelets due to
Daubechies are good candidates.  Our analysis so far is based on
matrix versions of (\relaxref{6.2}) of the form
\relaxeqb$$
C(z)\Psi (z^N)=\lambda\Psi(z),\quad z\in\bT,\ \lambda 
\mbox{ some constant matrix,}
\tag{6.4}
\relaxeqe$$
where $C$ is related to (\relaxref{6.1}), and $\Psi$ is a
matrix-valued function.  Thus as an added problem for future research,
we will study further the unitary part in the Wold decomposition of
$L^2(\bT)$ corresponding to the given isometric operator
\relaxeqb$$
\xi\mapsto {\frac 1{\sqrt{2}}} (m_0(z)\xi(z^2)+m_1(z)\xi(-z^2))\ .\tag{6.5}
\relaxeqe$$

Our preliminary examination indicates that the wavelets, which
correspond to the pairs $m_0,m_1$ (high pass/low pass filters) for
which the isometry in (\relaxref{6.5}) has a {nonzero unitary part} of
its Wold decomposition, are precisely the wavelets in $L^2(\bR)$ which
are equivalent to the familiar Haar wavelet.  But we plan to continue
and extend this research, as the idea appears to be also applicable
(with modifications and work) to other wavelets.  A second line of
research, connected with (\relaxref{6.5}), is to study the solutions
$\xi\ne 0$ which (for given $m_0,m_1$ as described above, see (6.1))
correspond to the unitary part of the Wold decomposition.  It turns
out that these solutions $\xi$ themselves generate {\em quadrature
mirror filters} and therefore correspond to orthogonal wavelets in
$L^2(\bR)$.  {We hope later to clarify this new form of duality for
wavelets in $L^2(\bR)$.}  As the idea seems basic, it should also be
useful (with further modifications) for understanding wavelets in
$L^2(\bR^d)$ when $d>1$.

\subsubhead
Acknowledgement
\endsubsubhead
Thanks are due to Reinhard Werner for permitting us to publish some results
from \cite{3}.

\end{document}